\documentclass[final,1p,times]{elsarticle}
\usepackage{amsmath,amsfonts,amssymb,amsthm}
\usepackage[all]{xy}
\usepackage{graphicx}
\usepackage{aliascnt,hyperref}

\journal{arXiv}

\newtheorem{teo}{Theorem}

\def\partet#1#2#3#4{\newaliascnt{#1}{#2}\newtheorem{#1}[#1]{#3}\aliascntresetthe{#1}\providecommand*{#4}{#3}}
\def\parted#1#2#3#4{\newaliascnt{#1}{#2}\newdefinition{#1}[#1]{#3}\aliascntresetthe{#1}\providecommand*{#4}{#3}}
\partet{lema}{teo}{Lemma}{\lemaautorefname}
\partet{prop}{teo}{Proposition}{\propautorefname}
\partet{cor}{teo}{Corollary}{\corautorefname}
\parted{defs}{teo}{Definition}{\defsautorefname}
\parted{ejemplo}{teo}{Example}{\ejemploautorefname}
\parted{obs}{teo}{Remark}{\obsautorefname}
\parted{parr}{teo}{}{\parrautorefname}
\parted{notc}{teo}{Notation}{\notcautorefname}
\newproof{dem}{Proof}

\def\C{\mathbb{C}}

\def\Q{\mathbb{Q}}
\def\PP{\mathbb{P}}

\def\then{\Longrightarrow}
\def\iff{\Longleftrightarrow}

\def\fl#1{\stackrel{#1}{\longrightarrow}}

\def\rk{\text{rk}}

\begin{document}

\begin{frontmatter}

\title{Quadratic equations of projective $PGL_2(\C)$-varieties.}

\author[dm]{C\'esar Massri\corref{correspondencia}\fnref{finanaciado}}
\address[dm]{Department of Mathematics, FCEN, University of Buenos Aires, Argentina}
\cortext[correspondencia]{Address for correspondence: Department of Mathematics, FCEN, University of Buenos Aires, Argentina}
\fntext[finanaciado]{The author was fully supported by CONICET, Argentina}
\ead{cmassri@dm.uba.ar}

\begin{abstract}
We make explicit the equations of any projective $PGL_2(\C)$-variety defined by quadrics.
We study their zero-locus and their relationship with the geometry of the Veronese curve.
\end{abstract}

\end{frontmatter}

\section*{Introduction.}

Due to the progress of mathematical computer systems, like Maple, Macaulay2, Singular, Bertini and others,
it is important to know explicitly the equations defining some known varieties. In this paper, we address this task
for projective varieties stable under $PGL_2(\C)$, the simplest of the simple Lie groups. In fact, we give
all the quadratic equations of any projective variety stable under $PGL_2(\C)$. We restrict ourselves to
varieties inside $\PP S^r(\C^2)$, where $r$ is a natural number.

Let $r\geq 2$ be a natural number.
Recall from \cite{MR1153249} that the $\mathfrak{sl}_2(\C)$-module $S^r(\C^2)$ is simple, that $S^r(\C^2)\cong S^r(\C^2)^\vee$
and that the decomposition of $S^2(S^r(\C^2))$ into simple submodules is given by
$$S^2(S^r(\C^2))=\bigoplus_{m\geq 0}S^{2r-4m}(\C^2).$$
In this article, we investigate varieties $M_m\subseteq\PP^r=\PP S^r(\C^2)$ generated in degree two by $S^{2r-4m}(\C^2)^\vee$.
Specifically, let $f_m:S^2(S^r(\C^2))\rightarrow S^{2r-4m}(\C^2)$ be the projection and let
$$M_m=\{x\in\PP S^r(\C^2)\,|\,f_m(xx)=0\}.$$
If $f_m=(q_0,\ldots,q_{2r-4m})$, then the generators of the ideal of $M_m$ are given by
$$\langle q_0,\ldots,q_{2r-4m}\rangle\cong S^{2r-4m}(\C^2)^\vee.$$

In the first section we study the equations defining $M_m$. In the second section we
give a bound for the dimension of the variety $M_m$. It is unknown if it is irreducible.
Any $PGL_2(\C)$-variety $X$ defined by quadrics is of the form
$$X=M_{m_1}\cap\ldots\cap M_{m_s},\quad I(X)_2=S^{2r-4m_1}(\C^2)^\vee\oplus\ldots\oplus S^{2r-4m_s}(\C^2)^\vee.$$
Then the knowledge of the quadratic equations of $M_m$ gives the explicit quadratic equations defining $X$.
Also, the bound on the dimension of $M_m$ gives a bound on the dimension of $X$.

{\bf Experiments and new theorems:} at the beginning of the first section, we found recursive equations that define
the quadrics containing the varieties $M_m$. We used a mathematical software (Maple) to compute
the coefficients of theses quadrics. All the observations made were proven in \autoref{carac-qk}, \autoref{simetria}
and \autoref{lema1}. We give a close formula for these coefficients.

Using the quadratic equations of the first section, we computed with another mathematical software (Macaulay2)
the dimensions and the degrees of $M_m$. We proved in \autoref{cota}, \autoref{tcr}, \autoref{hiper}, \autoref{complete} and \autoref{complete2}
some patterns that emerge from computations.

\section{Quadrics defining $M_m\subseteq\PP^r$.}
Let us fix a natural number $r$ and a projection $f_m:S^2(S^r(\C^2))\rightarrow S^{2r-4m}(\C^2)$.
For simplicity, let us denote $f=f_m$. Let $n=2r-4m$ be a fixed even number.

Consider the following basis in $\mathfrak{sl}_2(\C)$:
$$X=\left(\begin{array}{cc}0 & 1\\0 & 0\end{array}\right),\quad
H=\left(\begin{array}{cc}1 & 0\\0 & -1\end{array}\right),\quad
Y=\left(\begin{array}{cc}0 & 0\\1 & 0\end{array}\right).$$

Let $x_0\in S^r(\C^2)$ and $w_0\in S^n(\C^2)$ be maximal weight vectors.
The action of $Y\in\mathfrak{sl}_2(\C)$ on these vectors, generates bases
$\{x_0,\ldots,x_r\}$ of $S^r(\C^2)$ and $\{w_0,\ldots,w_n\}$ of $S^n(\C^2)$. Specifically,
$$x_i=\frac{Y^ix_0}{i!},\quad w_k=\frac{Y^kw_0}{k!},
\quad 0\leq i\leq r,\quad 0\leq k\leq n.$$
Using these bases, $f=\sum_0^n q_kw_k$, where $\{q_k\}_{k=0}^n$ are the quadratic equations of $M_m$.

Given that $f$ is $\mathfrak{sl}_2(\C)$-linear, we have the following relations:
$$Yf(x_ix_j)=f(Yx_ix_j)\iff\sum_{k=0}^nq_k(x_ix_j)Yw_k=
\sum_{k=0}^nq_k(Yx_ix_j)w_k\iff$$
$$\sum_{k=0}^{n-1}q_k(x_ix_j)(k+1)w_{k+1}=
\sum_{k=0}^nq_k((i+1)x_{i+1}x_j+(j+1)x_ix_{j+1})w_k\iff$$
$$kq_{k-1}(x_ix_j)=(i+1)q_k(x_{i+1}x_j)+(j+1)q_k(x_ix_{j+1}),\quad 0\leq k\leq n,\,0\leq i,j\leq r.$$
Note that all the forms depend recursively on $q_n$. In particular, if
$q_n=0$, the rest of the forms $q_k$ are zero.
Doing the same computation with $X$ instead of $Y$, we get a similar recursion:
$$(n-k)q_{k+1}(x_ix_j)=(r-i+1)q_k(x_{i-1}x_j)+(r-j+1)q_k(x_ix_{j-1}),\quad 0\leq k\leq n,\,0\leq i,j\leq r.$$
In these equations all the forms depend on $q_0$.
With $H$ we get conditions on each quadratic form,
$$Hf(x_ix_j)=f(Hx_ix_j)\iff\sum_{k=0}^nq_k(x_ix_j)Hw_k=\sum_{k=0}^nq_k(Hx_ix_j)\iff$$
$$\sum_{k=0}^nq_k(x_ix_j)(n-2k)w_k=\sum_{k=0}^nq_k((r-2i)x_ix_j+(r-2j)x_ix_j)w_k\iff$$
$$(n-2k)q_k(x_ix_j)=(2r-2(i+j))q_k(x_ix_j)\iff$$
$$(n-2k-2r+2i+2j)q_k(x_ix_j)=0,\quad 0\leq k\leq n,\,0\leq i,j\leq r.$$
Note that if $n-2r\neq 2k-2i-2j$, then $q_k(x_ix_j)=0$.
Saying thin in a different way, $q_k(x_ix_j)=0$ except maybe for $j=2m+k-i$.

\begin{cor}
Let $r$, $n$, $\{x_0,\ldots,x_r\}$ and $\{w_0,\ldots,w_n\}$ be as before and
let $q_0$ be an arbitrary bilinear form on $S^r(\C^2)$ such that:
$$0=(i+1)q_0(x_{i+1},x_j)+(j+1)q_0(x_i,x_{j+1}),\quad (2r-2i-2j-n)q_0(x_i,x_j)=0,\quad
0\leq i,j\leq r.$$
Then there exists a unique $\mathfrak{sl}_2(\C)$-morphism
$f:S^r(\C^2)\otimes S^r(\C^2)\rightarrow S^n(\C^2)$ such that its component over $w_0$ is $q_0$.
Even more, $f$ is symmetric if and only if $q_0$ is symmetric.
\end{cor}
\begin{dem}
Let $i,j,k$ be three integers such that $0\leq k\leq n,\,0\leq i,j\leq r$.
Assume we have defined $q_k$ and let us define $q_{k+1}$ using the recursive formula,
$$(n-k)q_{k+1}(x_i,x_j)=(r-i+1)q_k(x_{i-1},x_j)+(r-j+1)q_k(x_i,x_{j-1}).$$

Note that $q_{k+1}$ is symmetric if and only if $q_0$ is symmetric.
Let $f=q_0w_0+\ldots+q_nw_n$. By construction it is a $\mathfrak{sl}_2(\C)$-morphism and it is unique.\qed
\end{dem}

\begin{cor}\label{q_0}
A quadratic form $q_0$ that extends to an $\mathfrak{sl}_2(\C)$-map
$f:S^2(S^r(\C^2))\rightarrow S^{2r-4m}(\C^2)$, $f=q_0w_0+\ldots+q_nw_n$,
is given by
$$q_0(x_ix_{j})=
\begin{cases}
(-1)^i\binom{2m}{i}\lambda&\text{ if }j=2m-i\\
0&\text{otherwise}
\end{cases}$$
where $\lambda$ is a complex number.
In particular, if $\lambda\in\Q$, all the coefficients of $q_0$ are rational.
This implies that $q_k(x_ix_j)\in\Q$ for every $0\leq k\leq n$ and $0\leq i,j\leq r$.
\end{cor}
\begin{dem}
Let us analyze in more detail the hypothesis on the quadratic form $q_0$ given in the previous corollary.
The first condition,
$$0=(i+1)q_0(x_{i+1}x_j)+(j+1)q_0(x_ix_{j+1}),$$
implies that $q_0$ depends only on the values
$q_0(x_0x_j)$. This is because, given $q_0(x_0x_j)$ for every $0\leq j\leq r$, we may define
$$q_0(x_1x_j)=-\frac{j+1}{2}q_0(x_0x_{j+1}).$$
Thus, if we have defined up to $q_0(x_ix_j)$ for some $0<i<r$, we have
$$q_0(x_{i+1}x_j)=-\frac{j+1}{i+1}q_0(x_ix_{j+1}).$$

Let us discuss now the second hypothesis of the previous corollary,
$$(2r-2i-2j-n)q_0(x_ix_j)=0.$$
Given that $n=2r-4m$ we have $(2r-2i-2j-n)=0$ if and only if $i+j=2m$. Then
$$q_0(x_ix_j)\neq0\then i+j=2m.$$
Let $\lambda=q_0(x_0x_{2m})$ be arbitrary. Then applying the recursion we have
$$q_0(x_ix_{2m-i})=(-1)^i\binom{2m}{i}\lambda,\quad  0\leq i\leq 2m$$
\qed
\end{dem}

\begin{cor}
A $\mathfrak{sl}_2(\C)$-linear map $f:S^2(S^r(\C^2))\rightarrow S^{2r-4m}(\C^2)$ depends on one parameter, $\lambda\in\C$.
In other words,
$$\dim_{\mathfrak{sl}_2(\C)}(S^2(S^r(\C^2)),S^{2r-4m}(\C^2))=1.$$
\end{cor}
\begin{dem}
This fact is well known
but in this case we are emphasizing the fact that every morphism depends just on
one coefficient $\lambda$.\qed
\end{dem}

Now that we know exactly the coefficients of the quadratic form $q_0$,
let us study the other forms, $\{q_1,\ldots,q_n\}$.

First, we investigate the forms
$\{q_1,\ldots,q_{\frac{n}{2}}\}$. Then we prove that $q_k$ and $q_{n-k}$ are related ($0\leq k\leq\frac{n}{2}$).
\begin{teo}\label{carac-qk}
Let $\lambda=q_0(x_0x_{2m})$ and $j=2m+k-i$. Then for $0\leq k\leq\frac{n}{2}$,
$$\binom{n}{k}q_k(x_ix_j)=
\lambda\sum_{s=\max(0,i-k)}^{\min(2m,i)}(-1)^s\binom{2m}{s}\binom{r-s}{r-i}\binom{r-2m+s}{r-j}$$
\end{teo}
\begin{dem}
Recall these identities:
$$Xx_ix_j=(r-i+1)x_{i-1}x_j+(r-j+1)x_jx_{j-1}.$$
$$X^sx_i=(r-i+1)(r-i+2)\ldots(r-i+s)x_{i-s}=s!\binom{r-i+s}{r-i}x_{i-s}.$$
$$X^kx_ix_j=\sum_{l=0}^k\binom{k}{l}(X^lx_i)(X^{k-l}x_j).$$
From the equation $Xf(x_ix_j)=f(Xx_ix_j)$, we get
$$(n-k+1)q_k(x_ix_j)=q_{k-1}(Xx_ix_j).$$
Then
$$(n-k+1)(n-k+2)\ldots(n)q_k(x_ix_j)=(n-k+2)\ldots(n)q_{k-1}(Xx_ix_j)=$$
$$=(n-k+3)\ldots(n)q_{k-2}(X^2x_ix_j)=\ldots=q_0(X^kx_ix_j).$$

Without loss of generality we may assume $r>2m$.
When $r=2m$ (i.e. $n=0$) we obtain only $q_0$ that we already know (\autoref{q_0}).
Then
$$k!\binom{n}{k}q_k(x_ix_j)=q_0(X^kx_ix_j)=\sum_{l=0}^k\binom{k}{l}q_0(X^lx_iX^{k-l}x_j)=$$
$$=\sum_{l=0}^k\binom{k}{l}l!\binom{r-i+l}{r-i}(k-l)!\binom{r-j+k-l}{r-j}q_0(x_{i-l}x_{j-k+l})=$$
$$=\sum_{l=0}^k\binom{k}{l}l!\binom{r-i+l}{r-i}(k-l)!\binom{r-j+k-l}{r-j}(-1)^{i-l}\binom{2m}{i-l}\lambda.$$
Dividing by $k!$, the binomial $\binom{k}{l}$ simplifies.

Finally, making the change of variable $s=i-l$, we get
$$\binom{n}{k}q_k(x_ix_j)=\lambda\sum_{s=i-k}^i(-1)^s\binom{2m}{s}\binom{r-s}{r-i}\binom{r-2m+s}{r-j}.$$
By convention, the binomials that do not make sense are zero.\qed
\end{dem}

Let us prove now the relationship between the forms $q_k$ and $q_{n-k}$.
\begin{prop}\label{simetria}
Let $k$ and $i$ be two integers such that $0\leq k\leq r-2m$ and $0\leq i\leq r$. Let $j=2m+k-i$ and let $n=2r-4m$.
Then
$$q_k(x_ix_j)=q_{n-k}(x_{r-i}x_{r-j}).$$
\end{prop}
\begin{dem}
Recall the three conditions obtained from the fact that $f$ is $\mathfrak{sl}_2(\C)$-linear,
\begin{equation}\label{eq1}
kq_{k-1}(x_ix_j)=(i+1)q_k(x_{i+1}x_j)+(j+1)q_k(x_ix_{j+1}).
\end{equation}
\begin{equation}\label{eq2}
(n-k)q_{k+1}(x_ix_j)=(r-i+1)q_k(x_{i-1}x_j)+(r-j+1)q_k(x_ix_{j-1}).
\end{equation}
\begin{equation}\label{eq3}
(n-2k)q_k(x_ix_j)=(2r-2(i+j))q_k(x_ix_j).
\end{equation}
Let us make the following change of variables in the second recursion, (\ref{eq2}),
$$k'=n-k,\,i'=r-i,\,j'=r-j.$$
Note that $0\leq k'\leq n/2$ and $0\leq i',j'\leq r$. Then
\begin{equation}\tag{2'}\label{eq2p}
k'q_{k'-1}(x_{i'}x_{j'})=(i'+1)q_{k'}(x_{i'+1}x_{j'})+(j'+1)q_{k'}(x_{i'}x_{j'+1}).
\end{equation}

Let $a_k(i,j)=q_k(x_ix_j)$ and $b_{k'}(i',j')=q_{k'}(x_{i'}x_{j'})$. Then
\begin{equation*}\tag{1}
ka_{k-1}(i,j)=(i+1)a_k(i+1,j)+(j+1)a_k(i,j+1).
\end{equation*}
\begin{equation*}\tag{2'}
kb_{k-1}(i,j)=(i+1)b_k(i+1,j)+(j+1)b_k(i,j+1).
\end{equation*}
Then the recursions are the same.
If the initial data of these are equal, $a_{\frac{n}{2}}=b_{\frac{n}{2}}$,
then $q_k(x_ix_j)=q_{n-k}(x_{r-i}x_{r-j})$.
$$a_{\frac{n}{2}}(i,2m+\frac{n}{2}-i)=q_{\frac{n}{2}}(x_ix_{2m+\frac{n}{2}-i})=
q_{\frac{n}{2}}(x_ix_{2m+r-2m-i})=q_{\frac{n}{2}}(x_ix_{r-i})=$$
$$q_{\frac{n}{2}}(x_{r-i}x_i)=b_{\frac{n}{2}}(i,r-i)=b_{\frac{n}{2}}(i,2m+\frac{n}{2}-i).$$
\qed
\end{dem}

\begin{cor}
For every $0\leq k\leq n/2$ we have $\rk(q_k)=\rk(q_{n-k})\leq 2m+k+1$.
\end{cor}
\begin{dem}
The matrix assigned to the quadratic form $q_k$ has at least $2m+k+1$ nonzero coordinates.
They appear in some anti-diagonal ($i+j=2m+k$) making nonzero rows linearly independent.\qed
\end{dem}

In general, the equality does not hold. For example, if $r=6$ and
$n=4$ (that is, $m=2$), then $q_2(x_1x_5)=q_2(x_5x_1)=0$
making the rank less than or equal to $2+4+1$.
In this case, $\rk(q_0)=\rk(q_4)=5$, $\rk(q_1)=\rk(q_3)=6$ and $\rk(q_2)=5<7$.

Finally, let us give a lemma that we are going to use in the next section.
\begin{lema}\label{lema1}
Let $\lambda=q_0(x_0x_{2m})\neq 0$ and let $k$ be such that $0\leq k\leq n/2$. Then
$$q_k(x_0x_{2m+k})=q_{n-k}(x_{r}x_{r-2m})\neq 0.$$
Even more, if $m=0$,
$$q_k(x_ix_{k-i})=q_{n-k}(x_{r-i}x_{r-k+i})\neq 0,\quad 0\leq i\leq r.$$
\end{lema}
\begin{dem}
From \autoref{carac-qk} we have the formula
$$q_k(x_0x_{2m+k})=\lambda\frac{\binom{r-2m}{k}}{\binom{n}{k}}\neq 0.$$
And from \autoref{simetria},
$q_{n-k}(x_{r}x_{r-2m})=q_k(x_0x_{2m+k})\neq 0$.

Similarly if $m=0$,
$$q_{n-k}(x_{r-i}x_{r-k+i})=q_k(x_ix_{k-i})=\lambda\frac{\binom{r}{r-i}\binom{r}{r-k+i}}{\binom{n}{k}}\neq 0,\quad 0\leq i\leq r.$$\qed
\end{dem}

\section{Geometric properties of $M_m\subseteq\PP^r$.}
In the previous section we computed the equations for $M_m$. Recall that
$M_m\subseteq\PP S^r(\C^2)$ is a projective $PGL_2(\C)$-variety generated in degree two by
$$\langle q_0,\ldots,q_{2r-4m}\rangle\subseteq S^2(S^r(\C^2)^\vee).$$
In this section we use these equations to compute a bound for the dimension of $M_m$.

Let us introduce some new notation. Let
$$b_i^k(m)=b_i^k:=q_k(x_ix_{2m+k-i})=q_{n-k}(x_{r-i}x_{r-2m-k+i}),\quad 0\leq k\leq\frac{n}{2},\,0\leq i\leq r.$$
Given that $q_k$ is symmetric, we have $b_i^k=b_{2m+k-i}^k$.

If $x=a_0x_0+\ldots+a_rx_r$, then
$$q_k(a_0,\ldots,a_r)=\sum_{i=0}^{2m+k}q_k(x_ix_{2m+k-i})a_ia_{2m+k-i}=\sum_{i=0}^{2m+k}b_i^ka_ia_{2m+k-i}.$$
$$q_{n-k}(a_0,\ldots,a_r)=
\sum_{i=0}^{2m+k}q_{n-k}(x_{r-i}x_{r-2m-k+i})a_{r-i}a_{r-2m-k+i}=\sum_{i=0}^{2m+k}b_i^ka_{r-i}a_{r-2m-k+i}.$$
With this notation, let us write the derivatives of $q_k$ with respect to $a_i$,
$$\frac{\partial q_k(a_0,\ldots,a_r)}{\partial a_i}=b_i^ka_{2m+k-i}+b_{2m+k-i}^ka_{2m+k-i}=2b_i^ka_{2m+k-i}.$$

\begin{prop}\label{cota}
The variety $M_m\subseteq \PP^r$ has dimension $\dim(M_m)<2m$. If $m=0$, $M_m=\emptyset$.
\end{prop}
\begin{dem}
Let us compute the rank of the Jacobian matrix of
$$(a_0,\ldots,a_r)\rightarrow (q_0(a_0,\ldots,a_r),\ldots,q_n(a_0,\ldots,a_r)).$$
It is a $(n+1)\times(r+1)$-matrix.
{\footnotesize
$$\begin{pmatrix}
b_0^0a_{2m} & b_1^0a_{2m-1} & \ldots & b_{2m}^0a_0 & 0 & 0& 0& \ldots& 0 \\
b_0^1a_{2m+1} & \ldots & \ldots & \ldots & b_{2m+1}^1a_0 & 0 & 0&  \ldots & 0 \\
b_0^2a_{2m+2} & \ldots & \ldots & \ldots & \ldots & b_{2m+2}^2a_0 & 0 &  \ldots & 0 \\
& & & \vdots\\
b_0^{r-2m}a_{r} & b_1^{r-2m}a_{r-1} & \ldots & \ldots & \ldots & \ldots& \ldots &  \ldots & b_{r}^{r-2m}a_0\\
0 & b_0^{r-2m-1}a_r & b_1^{r-2m-1}a_{r-1} & \ldots & \ldots & \ldots& \ldots &  \ldots & b_{r-1}^{r-2m-1}a_1\\
0 & 0 & b_0^{r-2m-2}a_r & b_1^{r-2m-2}a_{r-1} & \ldots & \ldots& \ldots &  \ldots & b_{r-2}^{r-2m-2}a_2\\
& & & \vdots\\
0 & 0 & \ldots & \ldots & 0 & b_{0}^0a_{r} & b_{1}^0a_{r-1} & \ldots& b_{2m}^0a_{r-2m}\\
\end{pmatrix}$$
}

Let $Z$ be the hyperplane given by $\{a_r=0\}$.
From \autoref{lema1}, we know that $b_0^{k}\neq 0$ for $0\leq k\leq r-2m$. Then
for every point not in $Z$, the last $r-2m+1$
rows of the previous matrix are linearly independent
making the rank greater that or equal to $r-2m+1$. If $m=0$, the rank is $r+1$.

Take $X$ an irreducible component of $M_m$. It is also a $PGL_2(\C)$-variety.
Recall that the closure of an orbit must contain orbits of lesser dimension.
In particular, $X$ must contain a closed orbit.
The unique closed orbit of $PGL_2(\C)$ in $\PP S^r(\C^2)$ is the orbit of
the maximal weight vector, $x_0$, \cite[Claim 23.52]{MR1153249}.
Using the equivariant isomorphism $S^r(\C^2)\cong S^r(\C^2)^\vee$,
the vector $x_r$ corresponds to the maximal weight vector of $\PP S^r(\C^2)^\vee$. Then its orbit is closed in $\PP S^r(\C^2)^\vee$.
Applying the isomorphism again, we obtain a closed orbit in $\PP S^r(\C^2)$, hence
the orbit of $x_r$ is equal to the orbit of $x_0$.
This implies that the point corresponding to $x_r$, $(0:\ldots:0:1)$ is in $X$, hence $X\setminus Z$ is non-empty.
Then a generic smooth point of $X$ has dimension less than $2m$.
\qed
\end{dem}

\begin{notc}
Our intention now is to relate the geometry of the Veronese curve with the geometry of $M_m$.
This analysis gives a lower bound for the dimension of $M_m$.

Recall briefly the definition of the Veronese curve $c_r\subseteq\PP^r$ and its
osculating varieties $T^pc_r$.
The Veronese curve may be given parametrically (over an open affine subset) by
$$c_r:t\fl{}(1,t,t^2,\ldots,t^r).$$
Its tangential variety, denoted $T^1c_r$,
may be given by
$$(t,\lambda_1)\fl{}c_r+\lambda_1 c_r'.$$
It depends on two parameters. One indicates the point in the curve and the other, the tangent
vector on that point.

In general, its $p$-osculating variety, $T^pc_r$ is given by
$$(t,\lambda_1,\ldots,\lambda_p)\fl{}c_r+\lambda_1 c_r'+\ldots+\lambda_p c_r^{(p)}.$$
In each point of the curve, stands a $p$-dimensional plane.

We consider the curve $c_r$ and its osculating varieties $T^pc_r$ inside $\PP^r$.
The dimensions of $c_r$ and of $T^pc_r$ are the expected, $p+1$.

In the article \cite{MR994946}, the author computed the Hilbert polynomials of the varieties $T^p c_r$,
$$H_{T^p c_r}(d)=(dr-dp+1)\binom{p+d}{d}-(dr-dp+d-1)\binom{p+d-1}{d}.$$
This implies that $\dim(T^p c_r)=p+1$, $\deg(c_r)=r$ and $\deg(T^1 c_r)=2(r-1)$.
\end{notc}

\begin{prop}\label{tcr}
The variety $M_m$ contains $T^{m-1}c_r$ but does
not contain $T^mc_r$. In particular, $\dim(M_m)\geq m$.
\end{prop}
\begin{dem}
This proposition follows from \cite[Exercise 11.32]{MR1153249}.
It says that
$$I(T^{p}c_r)_2\cong\bigoplus_{\alpha\geq p+1}S^{2r-4\alpha}(\C^2).$$
Given that $S^{2r-4m}(\C^2)\subseteq I(T^{m-1}c_r)_2$ we get $I(M_m)\subseteq I(T^{m-1}c_r)$.

Similarly, if $I(M_m)_2\subseteq I(T^{m}c_r)_2$, then $S^{2r-4m}(\C^2)\subseteq I(T^{m}c_r)_2$. A contradiction.
\qed
\end{dem}

\begin{ejemplo}\label{hiper}
Suppose that $r$ is even and that $m=r/2$. Then we have exactly one equation $q_0$.
It is a quadratic form whose
matrix (diagonal of rank $r+1$) has coefficients $\lambda(-1)^i\binom{r}{i}$.
In fact this is the only quadric in $\PP^r$ invariant under $PGL_2(\C)$.
For $r=4$ this quadric is well known, \cite[10.12]{MR1182558}.

The variety $M_m=\PP\{q_0=0 \}\subseteq \PP^r$ is a quadric of maximal rank, hence irreducible.
Being a hypersurface, it has $\dim(M_m)=r-1$. Then, by \autoref{tcr}, we obtain
$$\begin{cases}
T^{\frac{r}{2}-1}c_r\subsetneq M_m & \text{if }r>2.\\
c_2=M_m & \text{if }r=2.
\end{cases}$$
With this example we deduce that the dimension of $M_m$ may be strictly greater than $m$.
\end{ejemplo}

\begin{teo}\label{complete}
If $r\geq 3$ is odd and $m=(r-1)/2$, then $M_m$ has codimension $3$ and degree $8$.
\end{teo}
\begin{dem}
We know that $I(M_m)=\langle q_0,q_1,q_2\rangle$ where
$$q_0(a_0,\ldots,a_r)=b_0^0a_0a_{r-1}+b_1^0a_1a_{r-2}+\ldots+b_{r-1}^0a_{r-1}a_{0},$$
$$q_1(a_0,\ldots,a_r)=b_0^1a_0a_{r}+b_1^1a_1a_{r-1}+\ldots+b_{r}^1a_{r}a_{0},$$
$$q_2(a_0,\ldots,a_r)=b_0^0a_{r}a_{1}+b_1^0a_{r-1}a_{2}+\ldots+b_{r-1}^0a_{1}a_{r}.$$
The coefficients of the quadratic forms satisfy the following relations
$$b_0^0=b_{r-1}^0,\quad b_1^0=b_{r-2}^0,\quad \ldots,\quad b_{m-1}^0=b_{m+1}^0,$$
$$b_0^1=b_{r}^1,\quad b_1^1=b_{r-1}^1,\quad \ldots,\quad b_{m-1}^1=b_{m+2}^1,\quad b_{m}^1=b_{m+1}^1.$$

To see that the dimension is $r-3$ let us compute the rank of the Jacobian matrix
at a specific point $p\in M_m$. The Jacobian matrix is given by
$$\begin{pmatrix}
b_0^0a_{r-1} & b_1^0a_{r-2} & \ldots & b_{r-1}^0a_0 & 0\\
b_0^{1}a_{r} & b_1^{1}a_{r-1} & \ldots &  b_{r-1}^{1}a_1 & b_{r}^{1}a_0\\
0 & b_{0}^0a_{r}  & \ldots& b_{r-2}^{0}a_{2}& b_{r-1}^{0}a_{1}\\
\end{pmatrix}.$$

Let $p=(p_0:\ldots:p_r)\in\PP^r$ be a point such that
$$p_i=\begin{cases}
1&\text{if }i=0\text{ or }i=m-1,\\
0&\text{otherwise.}
\end{cases}$$
Then $q_0(p)=q_1(p)=q_2(p)=0$, hence $p\in M_m$. The Jacobian matrix at $p$ is equal to
$$\left(\begin{array}{ccccccc}
0\ldots 0& b^0_{r-m}&    0        &    0      &0\ldots 0 & b_{r-1}^0&0 \\
0\ldots 0&    0      &  b^1_{r-m+1} &    0      &0\ldots0&   0          & b_{r}^{1}\\
0\ldots 0&    0      &     0       &b^0_{r-m+1} &0\ldots0&   0          & 0\\
\end{array}\right).$$
Given that $b_i^0\neq 0$ for all $0\leq i\leq r$ (see \autoref{q_0})
and that $b_r^1=b_0^1\neq0$ (see \autoref{lema1}) the previous matrix has maximal rank,
hence the codimension of $M_m$ at $p$ is equal to $3$. This implies that the codimension of $M_m$ is $3$ and the
degree is $8$.

Note that the point $p$ is in $T^{m-1}c_r$ and that the points on the curve $c_r$ are singular.\qed
\end{dem}

\begin{teo}\label{complete2}
If $r\geq 8$ is even and $m=r/2-1$, then $M_m$ has codimension $5$ and degree $32$.
\end{teo}
\begin{dem}
Let us argue as in the proof of \autoref{complete}. We know that $I(M_m)=\langle q_0,\ldots,q_4\rangle$,
$$q_0(a_0,\ldots,a_r)=\sum_{i=0}^{r-2} b^0_ia_ia_{r-2-i},\quad
q_1(a_0,\ldots,a_r)=\sum_{i=0}^{r-1} b^1_ia_ia_{r-1-i},$$
$$q_2(a_0,\ldots,a_r)=\sum_{i=0}^{r} b^2_ia_ia_{r-i},$$
$$q_3(a_0,\ldots,a_r)=\sum_{i=0}^{r-1} b^1_ia_{r-i}a_{i+1},\quad
q_4(a_0,\ldots,a_r)=\sum_{i=0}^{r-2} b^0_ia_{r-i}a_{i+2}.$$

Let $p=(p_0:\ldots:p_r)\in\PP^r$ be a point such that
$$p_i=\begin{cases}
1&\text{if }i=0\text{ or }i=m-1,\\
0&\text{otherwise.}
\end{cases}$$
Then $p\in M_m$. The Jacobian matrix at $p$ is equal to
$$\left(\begin{array}{cccccccccc}
0\ldots 0& b^0_{r-m-1} &     0        &    0      &0          &0          &0 \ldots0 & b_{r-2}^0    &     0    & 0\\
0\ldots 0&    0      &  b^1_{r-m} &    0      &0          &0          &0 \ldots0 &   0          & b_{r-1}^{1}& 0\\
0\ldots 0&    0      &  0           &b^2_{r-m+1}&0          &0          &0 \ldots0 &   0          & 0        &  b_{r}^{2}\\
0\ldots 0&    0      &     0        &    0      &b^1_{r-m+1}&0          &0 \ldots0 &   0          &    0&0\\
0\ldots 0&    0      &     0        &    0      &0          &b^0_{r-m+1}&0\ldots 0 &   0          &    0&0
\end{array}\right).$$
From \autoref{q_0} and \autoref{lema1}, we know that $b_0^2$, $b_0^1$, $b_{r-2}^0$
and $b^0_{r-m+1}$ are non-zero numbers.
But given that $b_r^2=b_0^2$ and $b_{r-1}^1=b_0^1$, they are also non-zero.
We need to prove that $b^1_{r-m+1}$ is non-zero for $r\geq 8$.
Recall that $b^1_{r-m+1}=b^1_{m-2}$.
$$b^1_{m-2}\neq 0\iff \binom{n}{1}q_1(x_{m-2}x_{m+3})\neq 0\iff
\sum_{s=m-3}^{m-2}(-1)^s\binom{2m}{s}\binom{r-s}{r-m+2}\binom{r-2m+s}{r-m-3}\neq 0\iff$$
$$\binom{2m}{m-3}(r-m+3)-\binom{2m}{m-2}(r-m-2)\neq 0\iff
\frac{m-2}{m+3}\neq\frac{r-m-2}{r-m+3}\iff$$
$$(m-2)(r-m+3)-(r-m-2)(m+3)\neq0\iff 10m-5r\neq0\iff 2m\neq r.$$
Given that $2m=r-2$, we obtain $b^1_{r-m+1}\neq 0$.\qed
\end{dem}

\begin{ejemplo}
We computed the dimension and the degree of $M_m$ for several values of $r$ and $m$:
$$\begin{array}{|c|c|c|c|c|c|c|c|c|c|c|c|c|}
\hline
m\backslash r & 2 & 3 & 4 & 5& 6& 7& 8& 9& 10 &11&12&13\\
\hline
1&\underline{1}&\underline{1}&\underline{1}&\underline{1}&\underline{1}&
\underline{1}&\underline{1}&\underline{1}&\underline{1}&\underline{1}&\underline{1}&\underline{1}\\
\hline
2 & & &\underline{3}&\underline{2}&3&2&2&2&2&2&3&2\\
\hline
3 &  & & & &\underline{5}&\underline{4}&\underline{3}&3&5&3&3&3\\
\hline
4 &  & & & & & &\underline{7}&\underline{6}&\underline{5}&4&4&4\\
\hline
5 &  & & & & & & & &\underline{9}&\underline{8}&\underline{7}&6\\
\hline
6 &  & & & & & & & & & &\underline{11}&\underline{10}\\
\hline
\end{array}$$
\begin{center}{\footnotesize Table: Dimension of $M_m\subseteq\PP^r$.}\end{center}

$$\begin{array}{|c|c|c|c|c|c|c|c|c|c|c|c|c|}
\hline
m\backslash r & 2 & 3 & 4 & 5& 6& 7& 8& 9& 10 &11&12&13\\
\hline
1&\underline{2}&\underline{3}&\underline{4}&\underline{5}&
\underline{6}&\underline{7}&\underline{8}&\underline{9}&\underline{10}&\underline{11}&\underline{12}&\underline{13}\\
\hline
2 & & &\underline{2}&\underline{8}&5&12&14&16&18&20&22&24\\
\hline
3 &  & & & &\underline{2}&\underline{8}&\underline{32}&21&12&27&30&33\\
\hline
4 &  & & & & & &\underline{2}&\underline{8}&\underline{32}&128&36&40\\
\hline
5 &  & & & & & & & &\underline{2}&\underline{8}&\underline{32}&128\\
\hline
6 &  & & & & & & & & & &\underline{2}&\underline{8}\\
\hline
\end{array}$$
\begin{center}{\footnotesize Table: Degree of $M_m\subseteq\PP^r$.}\end{center}

The numbers underlined are known in general (see \autoref{hiper}, \autoref{complete}, \autoref{complete2}).
Recall also that $m\leq\dim M_m<2m$.
\end{ejemplo}

\begin{obs}
To end this section, let us make a little remark and some more computations. Suppose now that
we want to study the variety $X$ defined by the quadrics that contain $T^pc_r$.
In other words, $X$ is generated in degree two and $I(X)_2=I(T^pc_r)_2$.

Given that $c_r$ is generated in degree two, when $p=0$, we have the equality, $X=c_r$.
In the general case, $T^pc_r\subseteq X$.

From \autoref{cota} and the fact that
$X=M_{p+1}\cap\ldots\cap M_{\lfloor r/2\rfloor}$, we get
$$p+1\leq\dim(X)\leq 2p+1.$$

We computed the dimension of the variety $X$ in the case $I(X)_2=I(T^pc_r)_2$:
$$\begin{array}{|c|c|c|c|c|c|c|c|c|c|c|c|}
\hline
&\PP^4&\PP^5&\PP^6& \PP^7& \PP^8& \PP^9& \PP^{10} &\PP^{11}&\PP^{12}&\PP^{13}\\
\hline
I(T^1c_r)_2&\underline{3}&\underline{2}&2&2&2&2&2&2&2&2\\
\hline
I(T^2c_r)_2& & &\underline{5}&\underline{4}&3&3&3&3&3&3\\
\hline
I(T^3c_r)_2& & & & &\underline{7}&\underline{6}&4&4&4&4\\
\hline
I(T^4c_r)_2& & & & & & &\underline{9}&\underline{8}&6&5\\
\hline
I(T^5c_r)_2& & & & & & & & &\underline{11}&\underline{10}\\
\hline
\end{array}$$
The dimensions underlined are those in which $I(T^pc_r)_2=I(M_m)_2$ for some $m$,
so it is information from a previous table.

In the variety $4$-osculating of $c_{12}\subseteq\PP^{12}$ the pattern breaks.
The dimension is $6$ instead of $5$. We deduce that this variety is not generated in degree two.

Assume now that $5\leq r\leq 8$. Let $X_r$ be the variety generated in degree two by $I(T^1c_r)_2$.
We computed that $X_r$ is irreducible, $\dim(X_r)=2$ and $\deg(X_r)=2(r-1)$.
Then we know explicitly the equations defining $T^1c_5$, $T^1c_6$, $T^1c_7$ and $T^1c_8$ (set-theoretically).

{\footnotesize
$$I(X_5)=\langle x_{5}x_{0}-3x_{4}x_{1}+2x_{3}x_{2},x_{4}x_{0}-4x_{3}x_{1}+3{x_{2}}^{2},x_{5}x_{1}-4x_{4}x_{2}+
3{x_{3}}^{2}\rangle.$$

$$I(X_6)=\langle x_{4}x_{0}-4x_{3}x_{1}+3{x_{2}}^{2},x_{6}x_{0}-9x_{4}x_{2}+8{x_{3}}^{2},x_{6}x_{2}-4x_{5}x_{3}+3{x_{4}}^{2},$$
$$x_{5}x_{0}-3x_{4}x_{1}+2x_{3}x_{2},x_{6}x_{1}-3x_{5}x_{2}+2x_{4}x_{3},x_{6}x_{0}-6x_{5}x_{1}+15x_{4}x_{2}-10{x_{3}}^{2}\rangle.$$

$$I(X_7)=\langle x_7x_3-4x_6x_4+3x_5^2,2x_7x_3+x_6x_4-3x_5^2,x_7x_2+3x_6x_3-4x_5x_4,x_3x_0-x_2x_1,$$
$$x_4x_0-4x_3x_1+3x_2^2,x_5x_0+3x_4x_1-4x_3x_2,x_7x_4-x_6x_5,2x_4x_0+x_3x_1-3x_2^2,$$
$$x_5x_0-3x_4x_1+2x_3x_2,x_6x_0-6x_5x_1+15x_4x_2-10x_3^2,x_6x_0-x_5x_1-5x_4x_2+5x_3^2,$$
$$x_6x_0+8x_5x_1+x_4x_2-10x_3^2,x_7x_0+5x_6x_1-21x_5x_2+15x_4x_3,x_7x_0+23x_6x_1+51x_5x_2-75x_4x_3,$$
$$x_7x_1+8x_6x_2+x_5x_3-10x_4^2,x_7x_1-x_6x_2-5x_5x_3+5x_4^2,x_7x_1-6x_6x_2+15x_5x_3-10x_4^2,$$
$$x_7x_2-3x_6x_3+2x_5x_4,x_7x_0-5x_6x_1+9x_5x_2-5x_4x_3,x_2x_0-x_1^2,x_7x_5-x_6^2\rangle.$$

$$I(X_8)=\langle x_4x_0-4x_3x_1+3x_2^2,x_8x_2-6x_7x_3+15x_6x_4-10x_5^2,x_8x_4-4x_7x_5+3x_6^2,$$
$$x_8x_1+2x_7x_2-12x_6x_3+9x_5x_4,x_8x_3-3x_7x_4+2x_6x_5,3x_6x_0-4x_5x_1-11x_4x_2+12x_3^2,$$
$$x_5x_0-3x_4x_1+2x_3x_2,x_7x_0+2x_6x_1-12x_5x_2+9x_4x_3,x_7x_0-5x_6x_1+9x_5x_2-5x_4x_3,$$
$$x_8x_1-5x_7x_2+9x_6x_3-5x_5x_4,x_6x_0-6x_5x_1+15x_4x_2-10x_3^2,$$
$$x_8x_0+12x_7x_1-22x_6x_2-36x_5x_3+45x_4^2,3x_8x_2-4x_7x_3-11x_6x_4+12x_5^2,$$
$$x_8x_0-2x_7x_1-8x_6x_2+34x_5x_3-25x_4^2,x_8x_0-8x_7x_1+28x_6x_2-56x_5x_3+35x_4^2\rangle.$$
}
\end{obs}

\section*{Acknowledgments.}
This work was supported by CONICET, Argentina.
The author thanks Fernando Cukierman, for useful discussions, ideas and suggestions.

\end{document}